\documentclass[12pt]{amsart}
\usepackage{amssymb}
\usepackage[all]{xy}

\theoremstyle{plain}
\newtheorem{theorem}{Theorem}[section]
\newtheorem{lemma}[theorem]{Lemma}
\newtheorem{corollary}[theorem]{Corollary}
\theoremstyle{definition}
\newtheorem{definition}[theorem]{Definition}
\theoremstyle{remark}
\newtheorem*{remark*}{Remark}

\newcommand{\gm}{\Gamma(Z)}
\newcommand{\hch}{\mathbf{H}}

\begin{document}
\title[A note on resolution]{A note on resolution of rational 
and hypersurface singularities}
\author{D.~A. Stepanov}
\address{The Department of Mathematical Modeling \\
Bauman Moscow State Technical University \\
Moscow 105005, Russia}
\email{dstepanov@bmstu.ru}
\thanks{The research was supported by RFBR, grant no. 05-01-00353, 
CRDF, grant no. RUM1-2692-MO-05, and the Program for the 
Development of Scientific Potential of the High School, 
no. 2.1.1.2381.}
\date{}

\begin{abstract}
It is well known that the exceptional set in a resolution of 
a rational surface singularity is a tree of rational curves.
We generalize the combinatoric part of this statement to higher
dimensions and show that the highest cohomologies of the dual
complex associated to a resolution of an isolated rational 
singularity vanish. We also prove that the dual complex associated 
to a resolution of an isolated hypersurface singularity is simply 
connected. As a consequence, we show that the dual complex
associated to a resolution of a 3-dimensional Gorenstein terminal
singularity has the homotopy type of a point.
\end{abstract}

\maketitle

\section{Introduction}
Let $o\in X$ be an isolated singularity of an algebraic variety
(or an analytic space) $X$ defined over a field of characteristic 0, 
$\dim X\geq 2$. Consider a good resolution $f\colon Y\to X$ (this 
means that the exceptional locus $Z\subset Y$ of $f$ is a divisor 
with simple normal crossings). Let $Z=\sum Z_i$, where $Z_i$ are
irreducible. To the divisor $Z$ we can associate the dual complex
$\gm$. It is a CW-complex whose cells are standard simplexes
$\Delta_{i_0\dots i_k}^{j}$ corresponding to the irreducible
components $Z_{i_0\dots i_k}^{j}$ of the intersections 
$Z_{i_0}\cap\dots\cap Z_{i_k}=\cup_j Z_{i_0\dots i_k}^{j}$. The
$k-1$-simplex $\Delta_{i_0\dots\widehat{i_s}\dots i_k}^{j'}$ is a
face of the $k$-simplex $\Delta_{i_0\dots i_k}^{j}$ iff 
$Z_{i_0\dots i_k}^{j}\cap Z_{i_0\dots\widehat{i_s}\dots i_k}^{j'}
\ne\varnothing$. If $X$ and $Y$ are surfaces, then $\gm$ is the
usual resolution graph of $f$. Note that $\gm$ is a simplicial
complex iff all the intersections $Z_{i_0}\cap\dots\cap Z_{i_k}$
are irreducible. This can be obtained for a suitable resolution.
Also note that if $\dim X=n$, then $\dim(\gm)\leq n-1$.

The complex $\gm$ was first studied by G.~L. Gordon in connection 
to the monodromy in families (see \cite{Gordon}). We say that $\gm$ 
is \emph{the dual complex associated to the resolution} $f$.
The main reason motivating the study of the dual complex is the fact 
that the homotopy type of $\gm$ depends only on the singularity
$o\in X$ but not on the choice of a resolution $f$. This is a
consequence of the Abramovich-Karu-Matsuki-W{\l}odarczyk Weak
Factorization Theorem in the Logarithmic Category (see \cite{AKMW}).
%and \cite{Matsuki}, Theorem~5-4-1).
Indeed, this theorem reduces the problem to the case of a blowup
$\sigma\colon(X',Z')\to(X,Z)$, where $X$ and $X'$ are smooth 
varieties with divisors $Z$ and $Z'$ with simple normal crossings,
and the center of the blowup is admissible in some sense. It 
can be explicitly verified that $\gm$ is homotopy equivalent to
$\Gamma(Z')$ (see \cite{Stepanov}).

For example (it is taken from \cite{Gordon}), consider the
singularity
$$\{x_{1}^{n+1}+x_{2}^{n+1}+\dots+x_{n}^{n+1}+x_1x_2\dots x_n=0\}
\subset\mathbb{C}^n\,.$$
A good resolution can be obtained just by blowing up the origin.
The reader can easily prove that the exceptional divisor $Z$ 
consists of $n$ hyperplanes in general position in 
$\mathbb{P}^{n-1}$. We see that the complex $\gm$ is the border
of a standard $n-1$-dimensional simplex and thus it has the 
homotopy type of the sphere $S^{n-2}$.

If $F\colon(Y,Z)\to(X,o)$ is a resolution of an isolated toric
singularity $(X,o)$, then the complex $\gm$ has the homotopy
type of a point (\cite{Stepanov}).

In this paper, we study the dual complex associated to a resolution
in the case when $X$ is a rational or a hypersurface singularity 
defined over the field $\mathbb{C}$ of complex numbers. We show that 
if $f\colon Y\to X$ is a good resolution of an isolated rational
singularity $o\in X$, $\dim X=n$, then $H^{n-1}(\gm,\mathbb{C})=0$
(see Theorem~\ref{T:rational}). The proof is a generalization
of M. Artin's argument from \cite{Artin} to the $n$-dimensional 
case. The main new ingredient is the lemma on the degeneracy of a
spectral sequence associated to a divisor with simple normal crossings
on a K{\"a}hler manifold (Lemma~\ref{L:ss}). If $X$ is an isolated
hypersurface singularity,  
$\dim X\geq 3$, then $\pi(\gm)=0$ (see Theorem 
\ref{T:hypersurface}). This result is based on the well known fact
that the link of an isolated hypersurface singularity of dimension
$n\geq 3$ is simply connected (\cite{Milnor}). These results
allow to prove that the homotopy type of the dual complex
associated to a resolution of an isolated rational hypersurface
3-dimensional singularity is trivial (Corollary~\ref{C:dim3}).
As an application, we show that the dual complex associated to a 
resolution of a 3-dimensional Gorenstein terminal singularity has 
the homotopy type of a point (Corollary~\ref{C:terminal}).

We prove our theorems for algebraic varieties but everything
holds also for analytic spaces (with obvious changes).

The author is grateful to V.~A. Iskovskikh, K. Matsuki, Yu. G.
Prokhorov, J. Steenbrink and J. Whal for written and oral 
consultations that were very stimulating and useful. Also we thank
to P. Popescu-Pampu who noticed a mistake in the first version
of the paper, to D. Arapura, P. Bakhtari, and J. W{\l}odarczyk who corrected
a mistake in the proof of Lemma~\ref{L:ss}, and to the referee of
\cite{Stepanov2} who suggested several improvements to the text.
% \pagebreak

\section{Rational singularities}\label{S:ratsing}
Recall the
\begin{definition}
An algebraic variety (or an analytic space) $X$ has \emph{rational 
singularities} if $X$ is normal and for any resolution 
$f\colon Y\to X$ all the sheaves $R^if_*\mathcal{O}_{Y}$ vanish, 
$i>0$.
\end{definition}

In the sequel, when we say that $f$ is a good resolution we 
additionally assume that all the intersections $Z_{i_0}\cap\dots
\cap Z_{i_k}$ of prime components of the exceptional divisor 
$Z=\sum Z_i$ of $f$ are irreducible, thus $\gm$ is a simplicial 
complex.

The following theorem can be considered as a generalization of 
the classical fact that the exceptional locus in a resolution of
a rational surface singularity is a tree of rational 
curves (\cite{Artin}).
\begin{theorem}\label{T:rational}
Let $o\in X$ be an isolated rational singularity of a variety 
(or an analytic space) $X$ of dimension $n\geq 2$, and let 
$f\colon Y\to X$ be a good resolution with the exceptional divisor 
$Z$. Then the highest (complex) cohomologies of the complex $\gm$ 
vanish:
% $$H^{n-1}(\gm,\mathbb{Z})=0\,.$$
$$H^{n-1}(\gm,\mathbb{C})=0\,.$$
\end{theorem}
\begin{proof}
Let $Z=\sum\limits_{i=1}^{N} Z_i$ be the decomposition of the 
divisor $Z$ to its prime components $Z_i$. We can assume that $X$ 
is projective (since the given singularity is isolated) and $f$
is obtained by a sequence of smooth blowups (Hironaka's resolution 
\cite{Hironaka}). Thus all $Z_i$ and $Y$ are K{\"a}hler manifolds.

The sheaves $R^i f_*\mathcal{O}_{Y}$  are concentrated at the 
point $o$. Via Grothendieck's theorem on formal functions 
(see \cite{Grothendieck}, (4.2.1), and \cite{Grauert}, Ch. 4,
Theorem~4.5 for the analytic case) the completion of the stalk of 
the sheaf $R^i f_*\mathcal{O}_{Y}$ at the point $o$ is
\begin{equation}\label{E:projlim}
\varprojlim_{(r)\to(\infty)} H^i(Z,\mathcal{O}_{Z_{(r)}})\,,
\end{equation}
where $(r)=(r_1,\dots,r_N)$ and $Z_{(r)}=\sum\limits_{i=1}^{N}
r_iZ_i$. If $(r)\geq (s)$, i.~e., $r_i\geq s_i$ $\forall i$, there 
is a natural surjective map $g$ of sheaves on $Z$: 
$$g\colon\mathcal{O}_{Z_{(r)}}\to\mathcal{O}_{Z_{(s)}}\,.$$ 
Since dimension of $Z$ is $n-1$, the map $g$ induces a surjective 
map of cohomologies
$$H^{n-1}(Z,\mathcal{O}_{Z_{(r)}})\to
H^{n-1}(Z,\mathcal{O}_{Z_{(s)}})\,.$$
Recall that the given singularity $o\in X$ is rational, and thus the
projective limit \eqref{E:projlim} is $0$. Therefore the 
cohomology group $H^{n-1}(Z,\mathcal{O}_Z)$ vanishes too (because
the projective system in \eqref{E:projlim} is surjective). Now it
follows from the lemma~\ref{L:cohomology} below that 
$H^{n-1}(\gm,\mathbb{C})=0$. 
% But since $\dim(\gm)\leq n-1$, we have
% also that $H^{n-1}(\gm,\mathbb{Z})=0$.
\end{proof}
\begin{lemma}\label{L:cohomology}
Let $Z=\sum Z_i$ be a reduced divisor with simple normal crossings
on a compact K{\"a}hler manifold $Y$, $\dim Y=n$, and assume that
$H^k(Z,\mathcal{O}_Z)=0$ for some $k$, $1\leq k\leq n-1$. Then the 
$k$-th cohomologies with coefficients in $\mathbb{C}$ of the complex
$\gm$ vanish too:
$$H^k(\gm,\mathbb{C})=0\,.$$
\end{lemma}
\begin{proof}
Let us introduce some notation. Let $Z^0=\sqcup_i Z_i$ be the 
disjoint union of the irreducible components $Z_i$ and $Z^p=
\sqcup_{i_0<i_1<\dots<i_p}Z_{i_0i_1\dots i_p}$ where 
$Z_{i_0i_1\dots i_p}=Z_{i_0}\cap Z_{i_1}\dots\cap Z_{i_p}$.
By $\varphi_p\colon Z^p\to Z$ denote the natural map. Consider the
structure sheaves $\mathcal{O}_{Z^p}=
\bigoplus\limits_{i_0<\dots<i_p}\mathcal{O}_{Z_{i_0\dots i_p}}$, 
the sheaves $\mathcal{A}^{p,q}=
\bigoplus\limits_{i_0<\dots<i_p}
\mathcal{A}^{p,q}_{Z_{i_0\dots i_p}}$
of differential forms of bidegree $(0,q)$ on $Z^p$ and their
direct images $\mathcal{K}^p=\varphi_{p*}\mathcal{O}_{Z^p}$ and
$\mathcal{K}^{p,q}=\varphi_{p*}\mathcal{A}^{p,q}$ on $Z$. The
sequence of sheaves $\{\mathcal{K}^p\}_{p\geq 1}$ forms a complex
via the combinatoric differentials $\delta^p\colon\mathcal{K}^p\to
\mathcal{K}^{p+1}$, where if 
$$a=\oplus a_{i_0\dots i_p}\in\mathcal{K}^p(U)=
\bigoplus\limits_{i_0<\dots<i_p}\varphi_{p*}
\mathcal{O}_{Z_{i_0\dots i_p}}(U)$$ 
is a section of the sheaf
$\mathcal{K}^p$ over an open set $U\subseteq Z$, then
$$(\delta(a))_{i_0\dots i_pi_{p+1}}(U)=\sum_{j=0}^{p+1}
(-1)^j(a_{i_0\dots\widehat{i_j}\dots i_{p+1}})
|_{Z_{i_0\dots i_{p+1}}\cap U}\,.$$
Note that there is also a natural injection of $\mathcal{O}_{Z}$
into $\mathcal{K}^0$. The sequence of sheaves
$$0\to\mathcal{O}_{Z}\to\mathcal{K}^0\overset{\delta^0}{\to}
\mathcal{K}^1\overset{\delta^1}{\to}\dots$$
is exact. This is easy to check by considering the stalks; 
in particular, the exactness at $\mathcal{K}^0$ is a consequence of
the following fact which holds locally in a sufficiently small 
neighborhood of every point of $Z$: if $\{f_i\}$ is a collection of 
regular functions on $Z_i$ such that their restrictions onto 
intersections $Z_i\cap Z_j$ coincide, then there exists a regular 
function $f$ on $Z$ such that $f|_{Z_i}=f_i$ for all $i$ (it is 
important here that the divisor $Z$ has normal crossings). Therefore 
the complexes
$$\mathcal{O}^*:\quad \mathcal{O}_{Z}\to 0\to 0\to\dots$$
and
$$\mathcal{K}^*:\quad \mathcal{K}^0\overset{\delta^0}{\to}
\mathcal{K}^1\overset{\delta^1}{\to}\mathcal{K}^2\overset{\delta^2}
{\to}\dots$$
are quasiisomorphic. It is clear that the hypercohomologies of
the first complex are isomorphic to the cohomologies of $Z$ with
coefficients in the structure sheaf: $\hch^p(\mathcal{O}^*)\simeq
H^p(Z,\mathcal{O}_{Z})$. Now let us calculate the hypercohomologies
of the complex $\mathcal{K}^*$ by using the acyclic 
resolutions
$$\mathcal{K}^p\to\mathcal{K}^{p,0}\overset{\bar{\partial}}{\to}
\mathcal{K}^{p,1}\overset{\bar{\partial}}{\to}\dots\,,$$
where $\bar{\partial}$ is the Dolbeaux differential.

Consider the bigraded sequence of groups $K^{p,q}=
H^0(\mathcal{K}^{p,q},Z)\simeq
\bigoplus\limits_{i_0<\dots <i_p}H^0(\varphi_{p*}\mathcal{A}^{p,q})$
endowed with the differentials $\bar{\partial}\colon K^{p,q}\to
K^{p,q+1}$ and $\delta\colon K^{p,q}\to K^{p+1,q}$, where
$\bar{\partial}$ is the Dolbeaux differential and $\delta$ is the
combinatoric differential defined as follows: if $\alpha=
\oplus\alpha_{i_0\dots i_p}\in K^{p,q}$, then
$$(\delta(\alpha))_{i_0\dots i_p i_{p+1}}=\sum_{j=0}^{p+1}
(-1)^{q+j}(\alpha_{i_0\dots\widehat{i_j}\dots i_{p+1}})
|_{Z_{i_0\dots i_{p+1}}}\,.$$
These differentials satisfy the equality $\bar{\partial}\delta+
\delta\bar{\partial}=0$, thus $(K^{*,*},\delta,\bar{\partial})$ is 
a bicompex. Let $(K^*,d)$, $K^n=\bigoplus\limits_{p+q=n}K^{p,q}$,
$d=\delta+\bar{\partial}$ be the associated complex. Now we can
write that $\hch^p(\mathcal{K}^*)=H^p(\mathcal{K}^*,d)$.

There is a filtration $F^p K^n=
\bigoplus\limits_{\substack{p'+q=n \\ p'\geq p}}K^{p',q}$ on the 
complex $K^*$. It is known (see \cite{GH}, Ch. 3, section 5) that 
the spectral sequence $E_r$ associated to the filtration $F^pK^n$ 
converges to the cohomologies $H^*(K^*,d)$ and
$$E_{0}^{p,q}\simeq K^{p,q}\,;$$
$$E_{1}^{p,q}\simeq H_{\bar{\partial}}^{q}(K^{p,*})\,;$$
$$E_{2}^{p,q}\simeq H_{\delta}^{p}(H_{\bar{\partial}}^{q}(K^{*,*}))
\,.$$
In particular, $E_{1}^{p,0}\simeq$$H_{\bar{\partial}}^{0}
(\bigoplus\limits_{i_0<\dots <i_p}H^0(\varphi_{p*}
\mathcal{A}^{p,*}))\simeq$$\bigoplus\limits_{i_0<\dots <i_p}
\mathbb{C}$. Therefore the cochain complex
$$0\to E_{1}^{0,0}\overset{\delta}{\to}E_{1}^{1,0}
\overset{\delta}{\to}\dots$$
is isomorphic to the cochain complex that one uses to 
calculate cohomologies of $\gm$ (here we denote by the same letter 
$\delta$ the map between cohomologies induced by the combinatoric
differential). It follows that
$$E_{2}^{p,0}\simeq H^{p}(\gm,\mathbb{C})\,.$$

We shall show that the spectral sequence $E_r$ degenerates in
$E_2$. The method of the proof is based on the standard technique of
the theory of the mixed Hodge structures. We learned it from
\cite{Kulikov}, Chapter 4, \S 2. Compare also \cite{Gordon2}. Since
the result about $E_r$ can be of a particular interest, we state it as
a separate lemma.
\begin{lemma}\label{L:ss}
Let $Z=\sum Z_i$ be a reduced divisor with simple normal crossings
on a compact K{\"a}hler manifold $Y$, and let $E_r$ be the associated
spectral sequence as described above. Then $d_r=0$ for all $r\geq 2$,
i.~e., this spectral sequence degenerates in $E_2$.
\end{lemma}
\begin{remark*}
Our proof of this lemma contained in \cite{Stepanov2} and in the previous
version of this preprint is incorrect. The mistake was that the restriction
of a harmonic differential form onto a subvariety is generally not harmonic 
even on compact K{\"a}hler varieties. The mistake was fixed by D. Arapura,
P. Bakhtary, and J. W{\l}odarczyk in \cite{ABW}, Corollary 2.3. In fact,
Lemma~\ref{L:ss} is a particular case of a more general result \cite{ABW},
Theorem 2.1. However, to keep our text self-contained we present a
corrected version of our original proof here.
\end{remark*}
\begin{proof}
Consider the diagram
$$
\xymatrix{
&  &  &  &  & \\
&\dots\ar[r]^\delta &K^{p,q}\ar[r]^\delta \ar[u]_{\bar{\partial}}
&K^{p+1,q}\ar[r]^\delta \ar[u]_{\bar{\partial}}
&K^{p+2,q}\ar[r]^\delta \ar[u]_{\bar{\partial}} &\dots \\
&\dots\ar[r]^\delta &K^{p,q-1}\ar[r]^\delta \ar[u]_{\bar{\partial}}
&K^{p+1,q-1}\ar[r]^\delta \ar[u]_{\bar{\partial}}
&K^{p+2,q-1}\ar[r]^\delta \ar[u]_{\bar{\partial}} &\dots \\
& &\ar[u]_{\bar{\partial}} &\ar[u]_{\bar{\partial}} 
&\ar[u]_{\bar{\partial}} &
}
$$
where $\bar{\partial}\delta+\delta\bar{\partial}=0$. First we
take cohomologies in the vertical rows and obtain the sequences
$$\dots\dots$$
$$
\dots\to H^q(K^{p,*})\overset{\delta}{\to} 
H^q(K^{p+1,*})\overset{\delta}{\to} 
H^q(K^{p+2,*})\to\dots 
$$
$$\dots\to H^{q-1}(K^{p,*})\overset{\delta}{\to} 
H^{q-1}(K^{p+1,*})\overset{\delta}{\to} 
H^{q-1}(K^{p+2,*})\to\dots
$$
$$\dots\dots$$
Here $\delta$ is the induced map between cohomologies. Then we
calculate $\delta$-cohomologies and obtain the differential
$$d_2\colon H_{\delta}^{p}(H_{\bar{\partial}}^{q}(K^{*,*})\to
H_{\delta}^{p+2}(H_{\bar{\partial}}^{q-1}(K^{*,*})$$
that acts as described below. 

An element $\bar{\bar{a}}\in H_{\delta}^{p}(H_{\bar\partial}^{q}
(K^{*,*}))$ is a class of those $\bar{a}\in 
H_{\bar\partial}^{q}(K^{p,*})$ modulo $\delta(H_{\bar\partial}^{q}
(K^{p-1,*}))$ that map to $0$ under the action of $\delta$: 
$\delta(\bar{a})=0$ in $H_{\bar\partial}^{q}(K^{p+1,*})$.
But $\bar{a}$ is a class of $a\in K^{p,q}\mod 
\bar\partial K^{p,q-1}$ such that $\bar\partial a=0\in K^{p,q+1}$.
Therefore we can choose a representative $a\in K^{p,q}$ for
$\bar{\bar{a}}$ such that $\bar\partial a=0$ and $\delta a=0$
modulo $\bar\partial K^{p+1,q-1}$ in $K^{p+1,q}$. It follows
that there exists an element $a'\in K^{p+1,q-1}$ such that
$\bar\partial a'=\delta a$. Map this $a'$ down to $K^{p+1,q-1}$:
$\delta(a')\in K^{p+2,q-1}$. It can be verified by standard methods 
that $\delta(a')$ determines correctly a class 
$\overline{\overline{\delta(a')}}$ in $H_{\delta}^{p+2}
(H_{\bar\partial}^{q-1}(K^{*,*}))$ and the differential $d_2$ is
defined as follows: 
$$d_2(\bar{\bar{a}})=\overline{\overline{\delta(a')}}\,.$$

Differentials $d_r$, $r\geq 3$, can be obtained by iterating this
construction. For example, $d_3$ is defined for $a\in K^{p,q}$ such
that $\delta(a')=0$ modulo $\bar\partial K^{p+2,q-1}$, thus there is
an $a''\in K^{p+2,q-2}$ such that $\bar\partial 
(a'')=\delta(a')$, and $d_3\colon E_{3}^{p,q}\to E_{3}^{p+3,q-2}$ is
induced by correspondence $a\to \delta(a'')$ (see, e.~g.,
\cite{Eisenbud}, A.3.13.4). 

Our aim is to show that
$d_r=0$, $r\geq 2$. The differential $d_2$ is trivial 
if the representative $a\in K^{p,q}$ can be chosen in such a way 
that $\delta(a)$ is exactly $0$ but not only $0$ modulo
$\bar\partial K^{p+1,q-1}$. But this is true because there are
harmonic differential forms in the class $\bar{a}\in 
H_{\bar\partial}^{q}(K^{p,*})$ and we can take $a$ to be a harmonic 
form of pure type $(0,q)$. Notice that since we work on a compact 
K{\"ahler} variety the form $a$ is not only $\bar{\partial}$-closed but 
also $\partial+\bar{\partial}$-closed, where $\partial+\bar{\partial}$ is 
the complex de Rham differential. The form $\delta{a}$ is defined by means 
of restrictions onto subvarieties and linear operations. We can not claim 
that it is also harmonic, but it remains $\partial+\bar{\partial}$-closed 
and of pure type $(0,q)$. But it is $0\mod(\bar\partial K^{p+1,q-1})$, 
i.~e.,  $\bar{\partial}$-exact. It follows from the 
$\partial\bar{\partial}$-lemma (\cite{GH}, Ch.~1, section~2) that 
$\delta{a}$ is also $\partial$-exact and thus is exactly $0$. Further, 
this $\delta a=0$ can be lifted to $a'=0$ in $K^{p+1,q-1}$, thus 
$\delta(a')=0$ and so forth. This shows that $d_r=0$ also for all $r\geq 3$.
\end{proof}

Now let us come back to the proof of Lemma~\ref{L:cohomology}. We have
$E_2=E_\infty$, therefore $H^{p}(\gm,\mathbb{C}) 
\simeq E_{2}^{p,0}$ is a subgroup of 
$$H^p(K^*,d)\simeq\hch^p(K^*)\simeq\hch^p(\mathcal{O}^*)\simeq
H^p(Z,\mathcal{O}_Z)\,.$$
If $H^p(\gm)$ is not trivial, then also $H^p(Z,\mathcal{O}_Z)$ is 
not trivial.
\end{proof}

\section{Hypersurface singularities}
\begin{theorem}\label{T:hypersurface}
Let $o\in X$ be an isolated hypersurface singularity of an
algebraic variety (or an analytic space) $X$ of dimension at least 
$3$ defined over the field $\mathbb{C}$ of complex numbers. If 
$f\colon Y\to X$ is a good resolution of $o\in X$, $Z$ its 
exceptional divisor, then the fundamental group of $\gm$ is trivial:
$$\pi(\gm)=0\,.$$
\end{theorem}
\begin{proof}
Let $n$ be dimension of $X$. We can assume that $X$ is a 
hypersurface in $\mathbb{C}^{n+1}$ and the singular point $o$
coinsides with the origin. Consider the link $M$ of singularity
$o\in X$, i.~e., the intersection of $X$ with a sphere $S^{2n+1}$
of sufficiently small radius around the origin. The link $M$ is
an $(n-2)$-connected smooth manifold (\cite{Milnor}, 
Corollary~2.9, Theorem~5.2), in particular, $M$ is simply connected.

We can also consider $M$ as the border of a tubular neighborhood
of the exceptional divisor $Z\subset Y$. It is known (see
\cite{ACampo}), that there is a surjective map $\varphi\colon
M\to Z$ whose fibers are tori. It follows that the induced map
$\varphi^*\colon\pi(M)\to\pi(Z)$ is also surjective and hence
$\pi(Z)=0$.

It remains to show that $\pi(Z)=0$ implies $\pi(\gm)=0$. It is
enough to construct a surjective map $\psi\colon Z\to\gm$
with connected fibers. The following lemma is, essancially, a partial
case of the general construction of a map from a topological space
$Z=\cup Z_i$ to the nerve $\gm$ corresponding to the covering
$\{Z_i\}$ (see \cite{Dold}, p. 355). As in 
section~\ref{S:ratsing}, we assume that the intersections 
$Z_{i_0\dots i_p}$ are irredusible so that $\gm$ is a simplicial 
complex.

\begin{lemma}
Let $Z$ be a divisor with simple normal crossings on an algebraic
veriety or an analytic space $X$, and let $\gm$ be the corresponding
dual complex. Then there exists a map $\psi\colon Z\to\gm$ which is
(i) simplicial in some triangulations of $Z$ and $\gm$, (ii)
surjective, and (iii) has connected fibers.
\end{lemma}
\begin{proof}
First, let us take a triangulation $\Sigma'$ of $Z$ such that all 
the intersections $Z_{i_0\dots i_p}$ are subcomplexes. Next we make 
the barycentric subdivision $\Sigma$ of $\Sigma'$ and the 
barycentric subdivision of the complex $\gm$. Now let $v$ be a vertex
of $\Sigma$ belonging to the  
subcomplex $Z_{i_0\dots i_p}$ but not to any smaller subcomplex 
$Z_{i_0\dots i_p i_{p+1}}$:
$$v\in Z_{i_0\dots i_p}\,,\quad v\notin Z_{i_0\dots i_p i_{p+1}}
\:\forall\,i_{p+1}\,.$$ 
Then let 
$$\psi(v)=\text{ the center of the simplex }
\Delta_{i_0\dots i_p}\,.$$
This determins the map $\psi$ completely
as a simplicial map (depending on the triangulation $\Sigma'$).
It is clear from our construction that $\psi$ is surjective.

We claim that the map $\psi$ has connected fibers. Indeed, first
observe that $\psi$ can be represented as a composition of 
topological contractions of connected subcomplexes
\begin{equation}\label{E:contr}
Z_{i_0\dots i_p}\setminus(\text{ the open neighborhoods
of intersections} 
\end{equation}
$$\text{of }Z_{i_0\dots i_p} \text{ with }Z_{i_{p+1}}\:
\forall i_{p+1}\ne i_0,\dots,i_p)\,.$$
By an open neighborhood of a subcomplex $K$ we mean the union of 
interior points of all simplicial stars of $\Sigma$ with centers
on $K$. All complexes in~\eqref{E:contr} are connected because
the codimension of $Z_{i_0\dots i_p}\cap Z_{i_{p+1}}$ in
$Z_{i_0\dots i_p}$ is $2$.

Further, the contraction of a subcomplex $K$ from \eqref{E:contr}
can be factored into one-by-one contraction of maximal
simplexes of $K$. The preimage of every connected set under such 
contraction is connected since the preimage of every simplex is a 
simplex. Therefore the map $\psi$ has all the needed properties.
\end{proof}
\end{proof}

Some important types of singularities are rational hypersurface.
Combining Theorems~\ref{T:rational} and \ref{T:hypersurface},
we can obtain some precise results in the $3$-dimensional case.
\begin{corollary}\label{C:dim3}
Let $o\in X$ be an isolated rational hypersurface singularity
of dimension $3$. If $f\colon Y\to X$ is a good resolution with the
exceptional divisor $Z$, then the dual complex $\gm$ associated
to the resolution $f$ has the homotopy type of a point.
\end{corollary}
\begin{proof}
We know from Theorems~\ref{T:rational} and \ref{T:hypersurface}
that $\gm$ is simply connected and $H^2(\gm,\mathbb{C})=0$. 
Since $\dim X=3$, we have $\dim(\gm)\leq 2$. Thus
$H_2(\gm,\mathbb{Z})=0$. Now Corollary~\ref{C:dim3} follows from 
the Inverse Hurevich and Whitehead Theorems.
\end{proof}

For instance, $3$-dimensional Gorenstein terminal singularities are 
exactly isolated compound Du Val points (up to an analytic 
equivalence, for details see \cite{Pagoda}). Here it is sufficient 
to us that compound Du Val points are hypersurface singularities.
On the other hand, terminal singularities (and, moreover, 
canonical) are rational (see \cite{Elkik}). Combining these
results with Corollary~\ref{C:dim3}, we get
\begin{corollary}\label{C:terminal}
Let $o\in X$ be a $3$-dimensional Gorenstein terminal singularity
and let $f\colon Y\to X$ be a good resolution with the 
exceptional divisor $Z$. Then the dual complex $\gm$ of $f$ has
the homotopy type of a point.
\end{corollary}

\end{document}